\documentclass{dimacs-l}

\title{Minimizing Polynomial Functions}
\author{Pablo A. Parrilo} 
\address{Department of Control and Dynamical Systems \\ 
California Institute of Technology\\
 Pasadena, CA 91125-8100.}
\email{pablo@cds.caltech.edu}

\author{Bernd Sturmfels}
\address{Department of Mathematics\\
 University of California\\
Berkeley, CA 94720.} 
\email{bernd@math.berkeley.edu}
\thanks{The collaboration between the authors was made possible
through a grant from the Center for Pure and Applied Mathematics at UC
Berkeley.
The second author was supported in part by NSF Grant 
DMS-9970254 and the Miller Institute at UC Berkeley}

\date{March 12, 2001.}


\usepackage{amsthm}
\usepackage{amsmath}
\usepackage{amsfonts} 
\newcommand{\baseRing}[1]{\ensuremath{\mathbb{#1}}}
\newcommand{\Z}{\baseRing{Z}}
\newcommand{\N}{\baseRing{N}}
\newcommand{\R}{\baseRing{R}}

\newcommand{\C}{\baseRing{C}} 

\theoremstyle{plain}
\newtheorem{theorem}{Theorem}[section]
\newtheorem{lemma}[theorem]{Lemma}
\newtheorem{corollary}[theorem]{Corollary}
\newtheorem{proposition}[theorem]{Proposition}

\newtheorem{remark}[theorem]{Remark}

\numberwithin{equation}{section}

\begin{document}


\begin{abstract} \noindent
We compare algorithms for global optimization of polynomial \linebreak
functions in many variables. It is demonstrated that existing
algebraic methods (Gr\"obner bases, resultants, homotopy
methods) are dramatically outperformed by a relaxation technique, due
to N.Z.~Shor and the first author, which involves sums of squares and
semidefinite programming. This opens up the possibility of using
semidefinite programming relaxations arising from  the Positivstellensatz
for a wide range of computational problems in real algebraic geometry.
\end{abstract}

\maketitle

\section{Introduction}

This is an expository and experimental paper concerned with the
following basic problem. Given a multivariate polynomial function $\,
f \in \R[x_1,\ldots,x_n]\,$ which is bounded below on $\R^n$, find the
global minimum $f^*$ and a point $p^*$ attaining it:
\begin{equation}
\label{DefProb}
 f^* \quad = \quad f(p^*) \quad = \quad
 {\rm min} \, \bigl\{ \, f(p) \, : \, p \in \R^n \bigr\}. 
\end{equation}
Exact algebraic algorithms for this task find all the critical points
and then identifying the smallest value of $f$ at any critical
point. Such methods will be discussed in Section 2. The techniques
include Gr\"obner bases, resultants, eigenvalues of companion matrices
\cite{clo}, and numerical homotopy methods \cite{li}, \cite{ver}.

An entirely different approach was introduced by N.Z.~Shor
(\cite{shor}, \cite{ss}) and further developed in the dissertation of
the first author \cite{pablo}.  The idea is to compute the largest
real number $\lambda$ such that $ f(x) - \lambda$ is a sum of squares in
$\R[x_1,\ldots,x_n]$.  Clearly, $\lambda$ is a lower bound for the optimal
value $f^*$.  We show in Section 3 that, when the degree of $f$ is
fixed, the lower bound $\lambda$ can be computed in polynomial time using
semidefinite programming \cite{VaB:96}.  If $\lambda = f^*$ holds then this
is certified by semidefinite programming duality, and the certificate
yields the optimal point $p^*$. In our computational experiments, to
be presented in Section 5, we found that $\lambda = f^*$ almost always
holds, and we solved  problems up to $n=15$.

The objective of this article is to provide a bridge between
mathematical programming and algebraic geometry, demonstrating
that algorithms from the former have the potential to play a major
role in future algorithms in the latter. This will be underlined in
Section 6, where we present open problems, and in Section 7
where we show that semidefinite programming in conjunction
with the Positivstellensatz is applicable to a wide range
of computational problems in real algebraic geometry.

\section{Computational Algebra}

In this section we discuss the following approach to our problem 
(\ref{DefProb}). Form the partial derivatives of
the given  polynomial  $f $ and consider the ideal they generate:
$$ I \quad = \quad 
\bigl\langle 
\frac{\partial f}{\partial x_1},
\frac{\partial f}{\partial x_2}, \ldots,
\frac{\partial f}{\partial x_n}
\bigr\rangle
\quad \subseteq \quad \R[x_1,\ldots,x_n] \, =: \, \R[{\bf x}]
$$
The zeros of the ideal $I$ in complex $n$-space
  $\C^n$ are the {\it critical points} of $f$.
Their number (counting multiplicity) is 
the dimension over $\R$ of the residue ring:
$$ \mu \quad = \quad
{\rm dim}_\R \,  \R[{\bf x}]/I  \quad = \quad \# \, \mathcal{ V}_\C (I). 
$$
We shall assume that $\mu$ is finite. (If $\mu = +\infty$ then
one can apply perturbation techniques to reduce to the 
case $\mu < + \infty$). For instance, if $f$ is a dense
polynomial of even degree $2d$ then it follows from
 {\it B\'ezout's Theorem} that
$\, \mu \, = \, (2d-1)^n $.

Consider the subset of {\it real critical points}:
$$ \mathcal{ V}_\R (I) \quad = \quad
\bigl\{ p^{(1)}, p^{(2)}, \ldots, p^{(\nu)} \bigr\} \quad \subset
\quad \R^n .$$
This set is  usually much smaller than the set of all complex critical points,
i.e.,  typically we have $\nu \ll \mu$. If we know the set
$ \mathcal{ V}_\R (I)$, then our problem is solved.

\begin{lemma}
\label{critical}
The optimal value is attained at a critical point:
\begin{equation}
f^* \quad = \quad
{\rm min} \bigl\{ \,
f(p^{(1)}),
f(p^{(2)}),\ldots,
f(p^{(\nu)}) \, \bigr\}
\end{equation}
\end{lemma}

The three techniques  to be described in this section
all compute  the set $ \mathcal{ V}_\R (I) $ of real critical points.
We will illustrate then for the following example:
\begin{equation}
\label{RunningEx}
   {\rm Minimize} \quad f(x,y,z) \,\,\, =   \,\,\,
 x^4  + y^4  + z^4  - 4 x y z + x + y + z
\end{equation}
The optimal value for this problem is $\, f^* \, = \, -2.112913882 $,
and, disregarding symmetry,
 there are three optimal points $(x^*,y^*,z^*)$ attaining this value:
$$            (0.988, -1.102, -1.102) \, , \,\,
            (-1.102, 0.988, -1.102) \, , \,\,
            (-1.102, -1.102, 0.988) .
$$

\subsection{Gr\"obner bases and eigenvalues}
\label{sec:eigs}
We review the method of solving polynomial equations by means of
Gr\"obner bases and eigenvalues \cite[\S 2.4]{clo}.  We are free to
choose an arbitrary term order $\prec$ on the polynomial ring
$\R[x_1,\ldots \! ,x_n]$. Let $\mathcal{ G}$ be a Gr\"obner basis for
the critical ideal $I$ with respect to $\prec$. While computing
Gr\"obner bases is a time-consuming task in general, this is not an
issue in this paper, since in all our examples the $n$ given
generators $\,{\partial f}/{\partial x_i}\,$ already form a Gr\"obner
basis in the total degree order. In our example the Gr\"obner basis is
\begin{equation}
\label{GrobnerBasis}
 \mathcal{ G} \quad = \quad \bigl\{ \,
\underline{x^3} - y z + 1/4 \, , \,\,
\underline{y^3} - x z + 1/4 \, , \,\,
\underline{z^3} - x y + 1/4 \,
\bigr\}
\end{equation}
A monomial $x_1^{u_1} \! \cdots x_n^{u_n}$ is {\it standard} if it is
not divisible by the leading term of any element in the Gr\"obner
basis $\mathcal{ G}$.  The set $\mathcal{ B}$ of standard monomials is
an $\R$-basis for the residue ring $\R[{\bf x}]/I$. The standard
monomials for (\ref{GrobnerBasis}) are:
\begin{equation}
\label{exampleB}
 \mathcal{ B} \quad = \quad \bigl\{ \,
x^i y^j z^k \, : \, i,j,k \in \N, \, i,j,k \leq 2  \,\bigr\}, \qquad
\mu \,=\, \# (\mathcal{ B})\, =\, 27. 
\end{equation}
For any polynomial $g \in \R[{\bf x}]$ consider the $\R$-linear 
endomorphism:
$$ {\rm Times}_g \, : \,\R[{\bf x}]/I \, \rightarrow \,  \R[{\bf x}]/I \, ,\,\,
h \,\mapsto \, g \cdot h.  $$
This endomorphism is represented in the basis $\mathcal{ B}$
by a real $\mu \times \mu$-matrix  ${\rm T}_g$. The entry
of ${\rm T}_g$ with row index ${\bf x}^u \in \mathcal{ B}$ and
column index ${\bf x}^v \in \mathcal{ B}$ is the coefficient
of ${\bf x}^v$ in the normal form of $\,{\bf x}^u \cdot g({\bf x})\,$
with respect to $\mathcal{ G}$.

\begin{proposition}
The optimal value $f^*$ is the smallest real eigenvalue of the matrix
${\rm T}_f$.  Any eigenvector of ${\rm T}_f$ with eigenvalue $f^*$
defines an optimal point $p^* = (p_1^*, \ldots, p_n^*)$ by the
eigenvector identities $\,{\rm T}_{x_i} \cdot v = p_i \cdot v \,$ for
$i=1,\ldots,n$.
\end{proposition}

\begin{proof}
This follows from  Lemma \ref{critical}
and Theorem (4.5) in the book of Cox-Little-O'Shea
\cite{clo}; see also
\cite[Exercise 17, page 62]{clo}.
\end{proof}

The resulting algorithm is to compute symbolically the matrices ${\rm
T}_f$ and ${\rm T}_{x_i}$ for $i=1,\ldots,n$, then compute numerically
its eigenvalues (and matching eigenvectors) of ${\rm T}_f$, and
finally determine $f^*$ and $p^*$ as in the proposition.

In our example the matrix ${\rm T}_f$ has format $27 \times 27$ with
rows and columns indexed by (\ref{exampleB}). Of its $729$ entries
only $178$ are nonzero. For instance, the column indexed by $xyz$ has
four nonzero entries, namely, the coefficients of
$$
{\rm normalform}_\mathcal{ G}( xyz \cdot f) \quad = \quad
\frac{3}{4} x^2 y z \, + \,
\frac{3}{4} x y^2 z \, + \,
\frac{3}{4} x y z^2 
\, - \, x^2 y^2 z^2. $$
The matrix ${\rm T}_f$ has maximal rank $27$. Of its eigenvalues only
three are real:
$$
-0.8692394998,\,
-0.8702981639, \,
-2.112913879 
$$
The three real eigenvalues have  multiplicity $3$,$ 1$, and $3$ respectively.

\subsection{Resultants and discriminants}

One algebraic method for solving polynomial equations is to use
resultants. Closely related to resultants are discriminants.
They express the condition on a hypersurface to have a singularity,
by means of a polynomial in the coefficients its defining equation.
 Let $t$ be a new indeterminate and form the discriminant of
the polynomial $\,f(x)-t\,$ with respect to $x_1,\ldots,x_n$:
$$
\delta(t) \quad := \quad
\Delta_x \bigl(\, f(x_1,\ldots,x_n) - t \, \bigr)
$$
Here $\Delta_x$ refers to the {\it $A$-discriminant}, defined in
\cite[Chapter 9]{gkzbook}, where $A$ is the support of $f$ together
with the origin. From \cite[\S 10.1.H]{gkzbook} we conclude that the
discriminant $\delta(t)$ equals the characteristic polynomial of the
matrix ${\rm T}_f$.

\begin{corollary}
The  optimal value $f^*$ is the smallest real root of  $\delta(t)$.
\end{corollary}

\smallskip

In our example,  $\delta(t)$ is
$256 t^3-512 t^2-96 t+473$ times the third power of
$$
65536 t^6+393216 t^5+1056768 t^4+1011712 
t^3-421376 t^2-437152 t+166419.
$$
The optimal value $f^* = -2.11...$ is a root of this sextic.
This sextic has Galois group $S_6$, so $f^*$ cannot be 
 expressed in radicals over the rationals.

The suggested algorithm is to compute $\delta(t)$, and minimal
polynomials for the coordinates $x_i^*$ of the optimal point, by
elimination of variables using matrix formulas for resultants and
discriminants \cite[Chapter 13]{gkzbook}.  The subsequent numerical
computation is to find the roots of a univariate polynomial.

\subsection{Homotopy methods}
\label{sec:hom}
The critical equations form a square system: $n$ equations in $n$
 variables having finitely many roots. Such a system is well-suited
 for {\it numerical homotopy continuation methods}.  For an
 introduction to this subject see the papers of Li \cite{li} and
 Verschelde \cite{ver}. The basic idea is to introduce a deformation
 parameter $\tau$ into the given system.  For instance, we might
 replace (\ref{GrobnerBasis}) by the following system which depends on
 a complex parameter $\tau$:
\begin{equation}
\label{homotopy}
H_\tau \,\,\, : \,\,\, 
x^3 - \tau \cdot y z + 1/4 \,= \,
y^3 - \tau \cdot x z + 1/4 \,= \,
z^3 - \tau \cdot x y + 1/4 \, = \, 0 .
\end{equation}
The solutions $\,\bigl( x (\tau),y (\tau),z (\tau) \bigr) \,$
are algebraic functions of $\tau$.
Our goal is to find them for $\tau = 1$.
It is easy to find the solutions for $\tau = 0$:
$$
(x(0),y(0),z(0))
\quad  = \quad \bigl(\,4^{-1/3} \cdot \eta_1,
\,4^{-1/3} \cdot \eta_2,\,
4^{-1/3} \cdot \eta_3 \,
\bigr) \, , \quad
\eta_i^3 = 1. 
$$
Homotopy methods trace the full set of solutions from
$\tau = 0$ to $\tau = 1$ along a suitable path in the
complex $\tau$-plane. 
We determine $f^*$ by evaluating the objective function
$f(x,y,z)$ at all the real solutions for $\tau = 1$.

Homotopy methods are frequently set up so that the system at $\tau =
0$ breaks up into several systems, each of which consists of
binomials.  If the input polynomials are sparse, then these are the
{\it polyhedral homotopies} which take the Newton polytopes of the
given equations into consideration. In the sparse case, the number
$\mu$ will be the {\it mixed volume} of the Newton polytopes.  For an
introduction to these polyhedral methods see \cite[Chapter 7]{clo} and
the references given there.

\subsection{How large is the B\'ezout number ?}

 A common feature of all three algebraic algorithms in this section is
that their running time is controlled by the number $\mu$ of complex
critical points.  In the eigenvalue method we must perform linear
algebra on matrices of size $\mu \times \mu$, in the discriminant
method we must find and solve a univariate polynomial of degree $\mu$,
and in the homotopy method, we are forced to trace $\mu$ paths from
$\tau = 0$ to $\tau = 1$. Each of these three methods becomes
infeasible if the number $\mu$ is too big; for instance,
$\,\mu \geq 10,000 \,$ might be too big.

Suppose that the given polynomial $f$ in $\R[x_1,\ldots,x_n]$ has
even degree $2d$ and is dense. This will be the case in the family of
examples studied in Section 5. Then $\mu$ coincides with the B\'ezout
number $ (2d-1)^n$.  Some small values for the B\'ezout number are
listed in Table \ref{tab:matsiz}.  Most entries in this table are
bigger than $10,000$.  We are led to believe that the algebraic
methods will be infeasible for quartics if $n \geq 8$.

Each entry in the first row of Table \ref{tab:matsiz} is a one.  This
means we can minimize quadratic polynomial functions by solving a
system of linear equations (in polynomial time).  The punchline of
this paper is to reduce our problem to a semidefinite programming
problem which can also be solved in polynomial time for fixed $d$.

\begin{table}
\begin{center}
\begin{tabular}{c|rrrrrrrr}
$2d \,\,\, \backslash \,\,\, n$ & 3 & 5 & 7 & 9 & 11 & 13  & 15 \\ \hline
2 &  1 &  1 &  1 &  1 &  1 &  1 &  1 \\
4 &  27 &  243 &  2187 &  19683 & 
177147 &
 1594323 &
... \\
6 &  125 &  3125 &  78125 &
 1953125 &
 48828125 & ... & ... \\
8 &  343 &  16807 &  
823543   & 
 40353607 & 
 1977326743 & ... & ... \\
10 &  729 &  59049 &  
 4782969 &
 387420489  & ... & ... & ... & \\
12 &  1331 & 
 161051 & 
19487171 &
 ... & ... & ... & ... \\
\end{tabular}
\end{center}
\caption{The B\'ezout number $\mu = (2d-1)^n$ for the critical equations.}
\label{tab:matsiz}
\end{table}

\section{Sums of Squares and Semidefinite Programming}
\label{sec:sossdp}

We present the method introduced by N.Z.~Shor (\cite{shor},
\cite{ss}), and further extended by the first author \cite{pablo}, for
minimizing polynomial functions.  This method is a relaxation: it
always produces a lower bound for the value of $f^*$. However, as we
shall see in Section 5, this bound very frequently agrees with $f^*$.

We may assume that the given polynomial $f(x_1,\ldots,x_n)$ has even
degree $2d$. Let $X$ denote the column vector whose entries are all
the monomials in $\,x_1,\ldots,x_n$ of degree at most $d$.  The length
of the vector $X$ equals the binomial coefficient
$$ N \,\, = \,\, \binom{n+d}{d} .$$
Let $\mathcal{ L}_f$ denote the set of all real symmetric
$N \times N$-matrices $A$ such that
$f({\bf x}) =  X^T \cdot A \cdot X$.
This is an affine subspace in the space of
real symmetric $N \times N$-matrices.
Assume that the constant monomial $1$ is the first entry of $X$.
Let $E_{11}$ denote the matrix unit whose only
nonzero entry is a one in the upper left corner.

\begin{lemma}
\label{lemma1}
For any real number $\lambda$, the following two are equivalent:
\begin{itemize}
\item The polynomial
$ \, f({\bf x}) - \lambda \,$ is a sum of squares in $\R[{\bf x}]$.
\item There is  a matrix $ A \in \mathcal{ L}_f$ such that
$\,A - \lambda \cdot E_{11} \,$ is positive semidefinite,
that is, all eigenvalues of this symmetric matrix
are non-negative reals.
\end{itemize}
\end{lemma}

\vskip -.1cm

\begin{proof}
The matrix $\,A - \lambda \cdot E_{11} \,$ is positive semidefinite
if and only if there exists a real Cholesky factorization
$\,  A - \lambda \cdot E_{11} \, = \, B^T \cdot B $. If
this holds then
$$
 f({\bf x}) - \lambda \, = \, 
  X^T \cdot A \cdot X - \lambda  \, = \, 
  X^T \cdot (  A - \lambda \cdot E_{11}) \cdot X \, = \, 
  X^T \cdot   B^T B \cdot X \, = \, 
(B X)^T \cdot (B X)
$$
is a sum of squares, and every sum of squares representation
arises in this way.
\end{proof}

We write $f^{sos}$ for the largest real number $\lambda$ for which the
two equivalent conditions are satisfied.  We always have $\, f^* \geq
f^{sos}$. This inequality may be strict. It is even possible that
$\,f^{sos} = - \infty $. An example of this form is Motzkin's
polynomial
\begin{equation}
m(x,y) \quad = \quad x^4 y^2 + x^2 y^4 - 3 x^2 y^2 .
\label{eq:motzkin}
\end{equation}
It satisfies $m(x,y) \geq -1$, but $m(x,y) - \lambda$
is not a sum of squares for \emph{any}  $\lambda \in \R$. We refer
to \cite{Reznick} for an excellent survey of the problem of
representing a polynomial as a sum of squares, and the important role
played by Motzkin's example.

Sums of squares are crucial for us because of the following complexity
result.

\begin{theorem}
\label{polytime}
Fix ${\rm deg}(f) = 2d$ and let the number of variables
$n$ vary. Then there exists a
polynomial-time algorithm, based on semidefinite programming, for
computing $f^{sos}$ from $f$. The same statement holds if $n$
is fixed and $d$ varies.
\end{theorem}

{\it Semidefinite programming} (SDP) is the study of optimization
problems over the cone of all positive semidefinite matrices.  This
branch of optimization has received a lot of attention in recent
years, both for its theoretical elegance and its practical
applications.  Semidefinite programs can be solved in \emph{polynomial
time}, using interior point methods; see \cite{NN}, \cite{HandSDP},
\cite{VaB:96}.  This complexity result (together with Lemma
\ref{lemma1}) implies Theorem \ref{polytime} because the quantity $\,N
= \binom{n+d}{d} = \binom{n+d}{n}\,$ grows polynomially if either $n$
or $d$ is fixed. This result appears in \cite{pablo}.

\begin{table}
\begin{center}
\begin{tabular}{c|rrrrrrrr}
$2d \,\,\, \backslash \,\,\, n$ & 3 & 5 & 7 & 9 & 11 & 13  & 15 \\ \hline
2 & 4 & 6 & 8 & 10 & 12 & 14 & 16 \\
4 &  10 & 21 & 36 & 55 & 78 & 105 & 136 \\
6 &  20 & 56 & 120 & 220 & 364 & 560 & 816  \\
8 &  35 & 126 &  330 & 715 & 1365 & 2380 & 3876 \\
10 &  56 & 252 & 792 & 2002 & 4368 & 8568 & 15504 \\
12 & 84 & 462 & 1716 & 5005 & 12376 & 27132 & 54264 \\
\end{tabular}
\end{center}
\caption{The matrix size $N = \binom{n+d}{d}$ for 
the semidefinite programs.}
\label{binocoeff}
\end{table}

Available implementations of interior-point methods for semidefinite
programming perform extremely well in practice, say, for problems
involving matrices up to $500$ rows and columns (provided there are not
too many variables). This allows for the efficient computation of
$f^{sos}$, and as we shall see in Section 4, SDP duality furnishes a
polynomial-time test to check whether $f^* = f^{sos}$ and for
computing the optimal point $p^*$ in the affirmative case.  A
comparison of Tables \ref{tab:matsiz} and \ref{binocoeff} suggests
that SDP has the potential to compute much larger instances than
algebraic methods. Section 5 will show that this is indeed the case.

\smallskip

Our example (\ref{RunningEx}) has parameters $d=2,n=3$. The
affine space $\mathcal{ L}_f$  consists of all
$10 \times 10$-matrices $A(\lambda,{\bf c})$
with $\lambda = 0$  and $c_i \in \R$ arbitrary in the family
{\smaller
$$
\left [\begin {array}{cccccccccc} -\lambda &
  1/2&  1/2&  1/2&  c_{1}&  -c_{2}&  
c_{3}&  -c_{4}&  -c_{5}&  c_{6}\\\noalign{\medskip}1/2&  -2\,c_{1}
&  c_{2}&  c_{4}&  0&  c_{7}&  -c_{8}&  -c_{9}&  c_{10}&  c_{12}
\\\noalign{\medskip}1/2&  c_{2}&  -2\,c_{3}&  c_{5}&  -c_{7}&  
c_{8}&  0&  c_{13}&  c_{14}&  c_{15}\\\noalign{\medskip}1/2&  c_{4}& 
c_{5}&  -2\,c_{6}&  c_{9}&  -c_{10}-c_{13}-2&  -c_{14}&  -
c_{12}&  -c_{15}&  0\\\noalign{\medskip}c_{1}&  0&  -c_{7}&  c_{9}&  1&  0&  
c_{16}&  0&  c_{17}&  c_{18}\\\noalign{\medskip}-c_{2}&  c_{7}&  
c_8 &  -c_{10}-c_{13}-2&  0&  -2\,c_{16}&  0&  -c_{17}&  -c_{19}&  
-c_{11}\\\noalign{\medskip}c_{3}&  -c_{8}&  0&  -c_{14}&  c_{16}
&  0&  1&  c_{19}&  0&  c_{20}\\\noalign{\medskip}-c_{4}&  -c_{9}&  
c_{13}&  -c_{12}&  0&  -c_{17}&  c_{19}&  -2\,c_{18}&  c_{11}&  0
\\\noalign{\medskip}-c_{5}&  c_{10}&  c_{14}&  -c_{15}&  c_{17}&  -
c_{19}&  0&  c_{11}&  -2\,c_{20}&  0\\\noalign{\medskip}c_{6}&  
c_{12}&  c_{15}&  0&  c_{18}&  -c_{11}&  c_{20}&  0&  0&  1\end {array}
\right ]
$$}
The rows and columns of this matrix  are indexed by 
the entries of the vector
$$
X \quad = \quad
\left [\begin {array}{cccccccccc} 1&x&y&z&{x}^{2}&xy&{y}^{2}&xz&yz&{z}
^{2}\end {array}\right ]^T.
$$
We invite the reader to check the identity
$$ X^T \cdot A(\lambda,{\bf c}) \cdot X \quad = \quad
f(x,y,z) \,-\, \lambda \qquad \text{for all} \,\,
c_1,\ldots,c_{20} \in \R.
$$
The lower bound $f^{sos}$ is the largest real number $\lambda $
such that, for some choice of $c_1,\ldots,c_{20} \in \R$, the
matrix $\,A({\bf c},\lambda) \,$ has all eigenvalues nonnegative.
We find that
$$ f^{sos} \quad = \quad f^* \quad = \quad  -2.112913882,$$
and the optimal matrix (to five digits) is given by:
{\tiny
$$
\left [\begin {array}{rrrrrrrrrr}
    2.1129 &   0.5000 &   0.5000 &   0.5000 &  -0.4678 &  -0.0922 &  -0.4678 &     -0.0922 &  -0.0922 &  -0.4678  \\
    0.5000 &   0.9356 &   0.0922 &   0.0922 &  -0.0000 &   0.0892 &  -0.0892 &	    0.0892 &  -0.6666 &  -0.0892  \\
    0.5000 &   0.0922 &   0.9357 &   0.0922 &  -0.0892 &   0.0892 &   0.0000 &	   -0.6667 &   0.0892 &  -0.0892  \\ 
    0.5000 &   0.0922 &   0.0922 &   0.9356 &  -0.0892 &  -0.6666 &  -0.0892 &	    0.0892 &   0.0892 &   0.0000  \\
   -0.4678 &  -0.0000 &  -0.0892 &  -0.0892 &   1.0000 &  -0.0000 &  -0.3180 &	    0.0000 &   0.0554 &  -0.3181  \\ 
   -0.0922 &   0.0892 &   0.0892 &  -0.6666 &  -0.0000 &   0.6360 &   0.0000 &	   -0.0554 &  -0.0554 &   0.0554  \\
   -0.4678 &  -0.0892 &   0.0000 &  -0.0892 &  -0.3180 &   0.0000 &   1.0000 &	    0.0554 &  -0.0000 &  -0.3180  \\
   -0.0922 &   0.0892 &  -0.6667 &   0.0892 &   0.0000 &  -0.0554 &   0.0554 &	    0.6361 &  -0.0554 &   0.0000  \\ 
   -0.0922 &  -0.6666 &   0.0892 &   0.0892 &   0.0554 &  -0.0554 &  -0.0000 &	   -0.0554 &   0.6360 &   0.0000  \\
   -0.4678 &  -0.0892 &  -0.0892 &   0.0000 &  -0.3181 &   0.0554 &  -0.3180 &	    0.0000 &   0.0000 &   1.0000  
\end{array} \right].
$$
} 
This matrix is positive semidefinite. By computing a factorization
$\,B^T \cdot B \,$ as in the proof of Lemma \ref{lemma1}, we can
express $f-f^{sos}$ as a sum of squares.  In the next section we show
how to recover the points at which the optimal value is achieved.

Note that the number $20$ of free parameters  is the
case ``$n=3, d=2$'' of:

\begin{remark}
\label{dimensionformula}
The dimension of  $\mathcal{L}_f$ equals the 
number of linearly independent quadratic relations
among the monomials of degree $\leq d$ in $n$ variables. It equals
\[
\mbox{dim } \mathcal{L}_f \quad = \quad \frac{1}{2} 
 \left [ \binom{n+d}{d} ^2 + \binom{n+d}{d}\right]  - \binom{n + 2d}{2d}.
\]
The codimension (with respect to the space of symmetric matrices) is
equal to
\[
\mbox{codim } \mathcal{L}_f \quad = \binom{n + 2d}{2d}.
\]

\end{remark}

\section{Semidefinite Programming Duality}

In Section 3 we demonstrated that computing $f^{sos}$ is equivalent to
minimizing a linear functional over the intersection of the affine
space $\mathcal{L}_f$ with the cone of positive semidefinite $N \times
N$-matrices.  In our discussion we have represented the space
$\mathcal{L}_f$ by a spanning set of matrices.  For numerical
efficiency reasons it is usually preferable to represent
$\mathcal{L}_f$ by its defining equations (unless $n$ and $d$ are very
small).

Duality is a crucial feature of semidefinite programming.  It plays an
important role in designing the most efficient interior-point
algorithms.  In what follows we review the textbook formulation of SDP
duality, in terms of matrices. Thereafter we present a reformulation
in algebraic geometry language, and we then explain how to test the
condition $f^{sos} = f^*$ and how to recover the optimal point $p^*$.

\subsection{Matrix Formulation}

Let $\mathcal{S}^N$ denote the real vector space
of symmetric $N\times N$-matrices, 
with the inner product $A \bullet B := \mbox{trace } (A
\,B)$, and the L\"owner partial order given by $A \preceq B$ if $B-A$
is \emph{positive semidefinite}. Recall that $A \in \mathcal{S}^N$ is
positive semidefinite if $x^T A x \geq 0$, for all $x \in \R^N$.
This condition is equivalent to nonnegativity of all 
eigenvalues of $A$, and to nonnegativity of all principal minors.

The general SDP problem  (\cite{VaB:96},
\cite{HandSDP}) can be expressed in the form:
\begin{equation}
\begin{array}{cccc}
\mbox{minimize } & F \bullet X &&\\
\mbox{subject to }& \mathcal{G} X &=& b \\
 &X & \succeq & 0
\end{array}
\label{eq:primalsdp}
\end{equation}
where $X,F \in \mathcal{S}^N$, $b \in \R^M$, and
$\mathcal{G}:\mathcal{S}^N \longrightarrow \R^M$ is a linear
operator. This is usually called the \emph{primal} form, in analogy
with the linear programming (LP) case.

Notice that (\ref{eq:primalsdp}) is a \emph{convex} optimization
problem, since the objective function is linear and the feasible set
is convex. There is an associated
\emph{dual} problem:
\begin{equation}
\begin{array}{cccc}
\mbox{maximize } &b^T y &&\\ \mbox{subject to }& F - \mathcal{G}^* y
\succeq 0
\end{array}
\label{eq:dualsdp}
\end{equation}
where $y \in \R^M$ and $\mathcal{G}^*:\R^M \longrightarrow
\mathcal{S}^N$ is the operator adjoint to $\mathcal{G}$. Any
feasible solution of the dual problem is a lower bound of the optimal
value of the primal:
$$
F \bullet X - b^T y \,\,\, = \,\,\,
F \bullet X - y^T \mathcal{G} X \,\,\, = \,\,\,
(F - \mathcal{G}^* y )\bullet X \geq 0.
$$
The last inequality holds since the inner product of two positive
semidefinite matrices is nonnegative.  The converse statement (primal
feasible solutions give upper bounds on the optimal dual value) is
obviously also true. The inequality above is called
\emph{weak duality}.  Under certain conditions (notably, the
existence of strictly feasible solutions), \emph{strong duality} also
holds: the optimal values of the primal and the dual problems
coincide. If strong duality holds, then at optimality the matrix
$\, X \cdot (F - \mathcal{G}^* y)\,$ is zero, since $A,B \succeq 0,
\mbox{trace}(AB)=0 $ implies $AB=0$. This can be interpreted as a
generalization of the usual \emph{complementary slackness} LP
conditions.

Practical implementations of SDP (we will use SeDuMi \cite{sedumi})
simultaneously compute both the optimal matrix $X$ for
(\ref{eq:primalsdp}) and the optimal vector $y$ for
(\ref{eq:dualsdp}).

\subsection{Polynomial Formulation}
We set $m = N-1$ and we identify $\mathcal{S}^N$ with the
real vector space $\,\R[x]_2 = \R[x_0,x_1,\ldots,x_m]_2 \,$
of quadratic forms in $m+1$ variables. The vector space
dual to $\mathcal{S}^N$ is now denoted
$\,\R[\partial]_2 =\R[\partial_0,\partial_1,\ldots,\partial_m]_2 $.
The dual pairing is  given by differentiation and is denoted $\bullet$.
For any $f \in \,\R[\partial]_2 \,$ and any 
 real vector $p = (p_0,\ldots,p_m) \in \R^{m+1}$, the following
familiar identity holds:
\begin{equation}
\label{pointEval}
f(\partial_0,\ldots,\partial_m) \bullet  \frac{1}{2} (\sum_{i=0}^m p_i x_i)^2
\quad = \quad f(p_0,\ldots,p_m).
\end{equation}
We consider the general quadratic programming problem:
\begin{equation}
\label{QP}
\hbox{Minimize} \quad f(p)  \quad
\hbox{subject to}\quad
g_0(p) = 1 
\,\, \hbox{and} \,\,
g_1(p) = \cdots =  g_r(p) =  0 ,
\end{equation}
where $f,g_0,\ldots,g_r \in \R[\partial]_2 $ are given
and we are looking for an  optimal point $p \in \R^{m+1}$.
This problem can be relaxed to the following 
{\bf primal SDP:}
\begin{eqnarray*}
& \hbox{Minimize} \quad f(\partial) \bullet q(x) \quad
 \hbox{subject to} \quad q(x) \succeq 0
\\
& \hbox{and} \quad
g_0(\partial) \bullet q(x) = 1
\,\, \hbox{and} \,\,
g_1(\partial) \bullet q(x) = \cdots = 
g_r(\partial) \bullet q(x) = 0.
\end{eqnarray*}
The inequality $\, q(x) \succeq 0 \,$ means that $q$ is non-negative
on $\R^{m+1}$, i.e., $q$ is in the positive semidefinite cone in
$\R[x]_2$. In view of (\ref{pointEval}), the optimal value of
(\ref{QP}) is greater than or equal to the optimal value of the primal
SDP, and equality holds if and only if there is an optimal solution of
the form $q(x) = \frac{1}{2} (\sum_{i=0}^m p_i x_i)^2$.

Every semidefinite programming problem comes with
a dual problem, as in the previous subsection; see also
\cite{VaB:96}. In our case the
{\bf dual SDP} takes the form:
\begin{eqnarray*}
& \hbox{Maximize the first coordinate } \,\,\, \lambda  \,\,\,
\hbox{of the vectors} \quad (\lambda, \mu_1,\ldots,\mu_r ) \in \R^{r+1} \\
& \hbox{subject to the conditions}
\quad
f(\partial) + \sum_{i=1}^r \mu_i \cdot g_i(\partial) - 
\lambda \cdot g_0(\partial) \,\succeq \, 0
\end{eqnarray*}
Assuming the existence of a strictly feasible primal solution, 
the maximum value in the {\bf dual SDP} is always equal to
the  minimum value in the {\bf primal SDP}. Under this 
regularity assumption, 
which is easy to satisfy in our application, we conclude:

\begin{proposition}
\label{completesquare}
If the primal SDP has an optimal solution of the form
$\,q(x) = \frac{1}{2} (\sum_{i=0}^m p_i x_i)^2\,$
then the vector $(p_0,\ldots,p_m)$ is an optimal
solution for (\ref{QP}).
\end{proposition}

\subsection{Minimizing Quadratic Functions over Toric Varieties}

A {\it toric variety} is an algebraic variety, in affine space 
or projective space, which has a parametric representation
by monomials. Equivalently, a toric variety is an irreducible
variety which is cut out by {\it binomial equations}, that is,
differences of monomials. Here we will be interested in
those projective toric varieties which are defined by
quadratic binomials. This class includes many examples
from classical algebraic geometry, such as
Veronese and Segre varieties.
See \cite{Stbook} for an introduction.

Let $X$ be a toric variety in projective $m$-space 
whose defining prime ideal 
is generated by quadratic binomials $g_1,\ldots,g_r$
in $\R[\partial_0,\ldots,\partial_m]$.
Each generator has the form 
$\, \partial_i \partial_j -  \partial_k \partial_l \,$
for some $i,j,k,l \in \{0,1,\ldots,m\}$.
We set $\,g_0 (\partial) = \partial_0^2 $.
Then the equation  $\,g_0 (\partial) = 1\,$ on $X$
defines an affine toric variety $\tilde X$, such that
$X$ is the projective closure of $\tilde X$.
Every quadratic polynomial function on the affine variety
$\tilde X$ is represented  by a quadratic form
$f \in \R[\partial]_2$ as above. This representation
is unique modulo the $\R$-linear span of $g_1,\ldots,g_r$.
Our problem (\ref{QP}) is hence equivalent to
minimizing a quadratic function
over an affine toric variety
defined by quadrics:

\begin{equation}
\label{toricQP}
\hbox{Minimize} \quad f(p)  \quad
\hbox{subject to}\quad p \in \tilde X 
\end{equation}

The optimal value  of the
{\bf dual SDP} relaxation in Subsection 4.2
is the largest real number $\lambda$ such that
$ f - \lambda$ is a sum of squares in the
coordinate ring of $\tilde X$.

Let us now return to our original problem (\ref{DefProb}) where the
given polynomial is dense of degree $2d$ in $n$ variables.  Here $X$
is the {\it Veronese variety} in projective $N$-dimensional space
which is parameterized by all monomials of degree at most $d$.  (If
the polynomial in (\ref{DefProb}) is sparse then another toric variety
can be used.) Writing our given polynomial as a quadratic form in 
homogeneous coordinates on $X$, our minimization problem
(\ref{DefProb}) is precisely the quadratic toric problem
(\ref{toricQP}).

We solve  (\ref{toricQP}) by simultaneously solving 
the primal and dual SDP relaxation in Subsection 4.2.
If the optimal value $\lambda$ of the dual SDP agrees with the true
minimum of $f$ over $\tilde X$ then the primal SDP
has an optimal solution
$\,q(x) = \frac{1}{2} (\sum_{i=0}^m p_i x_i)^2$
which exhibits an optimal point $(p_0,\ldots,p_m) \in X$
at which $f$ is minimized. 

In our running example, we have $m = 9$ and $r=20$, 
and $X$ is the 
quadratic Veronese three-fold in projective $9$-space
which is given parametrically as
$$ 
(x_0 : x_1 : \cdots : x_9) \quad = \quad
\bigl( 1 : r : s : t : r^2 : rs :  s^2 : rt : st : t^2 \bigr).
$$
It is cut out by twenty quadratic binomials such as
$\,x_0 x_5 - x_1 x_2 $. These binomials correspond to the
parameters $c_i$ in the $10 \times 10$-matrix $A(\lambda,{\bf c})$
in Section 2.

\section{Experimental Results}

We now present our computational experience with 
Shor's relaxation for global minimization of polynomial functions.
As mentioned earlier, the computational advantages of our method
are based on the following three independent facts:

\begin{itemize}
\item The dimension $N$ of the matrix required in the sum of squares
formulation is much smaller than the B\'ezout number $\mu$, since it only
scales polynomially with the number of variables.
See Tables \ref{tab:matsiz} and \ref{binocoeff} above.

\item Semidefinite programming provides an efficient 
algorithm for deciding whether a polynomial
is a sum of squares, and to find such representations
for polynomials whose coefficients may depend
linearly on  parameters.

\item The lower bound $f^{sos}$ very often coincides
with the  exact solution $f^*$ of our problem (\ref{DefProb}),
at least for the class of problems analyzed here.
\end{itemize}
The experimental results in this section strongly support the 
validity of these facts.

\subsection{The test problems}
For our computations, we fix a positive integer $K$,
and we sample from the following family of
polynomials of degree $2d$ in $n$ variables:
\begin{equation}
f(x_1,\ldots,x_n) \quad =  \quad
x_1^{2d} + x_2^{2d} + \cdots + x_n^{2d} \,+\, g(x_1,\ldots,x_n)
\label{eq:polyfam}
\end{equation}
where $g \in \Z[x_1,\ldots,x_n]$ is a random polynomial of total degree
$\leq 2d-1$ whose $\,\binom{n + 2d-1}{n}\,$
coefficients are independently and
uniformly distributed among integers  between $-K$ and $K$.
Thus our family depends on three parameters: $n$, $d$ and $K$.

This family has been selected to ensure three important properties:
\begin{description}
\item [Boundedness] The highest order terms 
$\, x_i^{2d} \,$ ensure that $f$
is bounded below, and that the minimum value $f^*$ is achieved at some
point $\,p^* \in \R^n$.
\item [Efficient basis computation]
When solving polynomial systems, the calculation
of a Gr\"obner basis is a time-consuming task.
The structure of the polynomial (\ref{eq:polyfam})
allows us to bypass this expensive step, since the set of $n$
scaled partial derivatives $\, x_i^{2d-1} + \frac{1}{2d} \cdot
\partial g / \partial x_i\,$ is already a  Gr\"obner basis
with respect to total degree; cf.~\cite[\S 2.9, Proposition 4]{clobaby}.
\item [Simplicity]
A main reason for this choice of model is its simplicity. 
While more sophisticated
choices have other desirable mathematical properties (such as
invariance under certain transformations), we preferred to
analyze here, as a first step, a relatively easy to describe set of
instances.
\end{description}

An important question is if the structure of the polynomials
(\ref{eq:polyfam}) is
somehow ``biased'' towards the application of sum of squares
methods. This is a relevant issue, since the performance
of algorithms on ``random instances'' sometimes
provides more information on the
problem family, rather than on the algorithm itself. 
Concerning this question, we limit ourselves to notice that, for $K$ 
sufficiently large, the family
(\ref{eq:polyfam})  does include polynomials $f$ with $f^{sos} < f^*$.
A simple example is $f(x,y) =
x^8 + y^8 + 2700 \, m(x,y)$, where $m(x,y)$ is the Motzkin
polynomial (\ref{eq:motzkin}).

The polynomials in our family have global minima that
generally have large negative values, of the order of $-K^{2d}$.  
This leads to ill-conditioning of the symmetric matrices
described in  Lemma \ref{lemma1},
and hence to numerical problems for the interior-point
algorithm. Our remedy is  a simple homogeneous scaling of the form
\[
f_s (x_1,\ldots,x_n) \quad  =  \quad
 \alpha^{-2d} \cdot f(\alpha x_1, \ldots, \alpha x_n) ,
\quad \qquad \hbox{for some $\alpha >0$}.
\]
Obviously, this does not affect the properties of being a sum of
squares, or whether $f^*=f^{sos}$. However, as is generally the rule
in numerical optimization, this scaling step greatly affects both the
speed and the accuracy of the SDP solution.

\subsection{Algorithms and software}

Most of the test examples were run on a Pentium III 733Mhz with 256
MB, running Linux version 2.2.16-3, and using MATLAB version
5.3. Because of physical memory limitations, our largest examples
(quartics in fifteen variables), were run on a Pentium III 650Mhz with
320 MB, under Windows 2000. The semidefinite programs were solved
using the SDP solver SeDuMi \cite{sedumi}, written by Jos Sturm. It is
currently one of the most efficient codes available, at least for the
restricted class of problems relevant here. SeDuMi can be run from
within MATLAB, and implements a self-dual embedding technique. The
default parameters are used, and the solutions computed are typically
exact to machine precision (SeDuMi provides an estimate of the quality
of the solution).

The MATLAB Optimization toolbox was used for the implementation of a
local search approach, to be described in Section~\ref{sec:local}. For
the numerical homotopy method, we used the software PHCpack
\cite{PHC}, written by Jan Verschelde. The computation of the sparse
matrix ${\rm T}_f$ was done using Macaulay 2 \cite{M2}, and its
eigenvalues were numerically computed using MATLAB.

We do not make strong claims about the efficiency of our
implementations: while reasonable fast, for large scale problems
considerable speedups are possible at the expense of customized
algorithms. Nevertheless, we believe that the  issues
raised regarding the applicability of algebra-based techniques to
problems with large B\'ezout number remain valid, independently of the
particular software employed.

\subsection{Standard local optimization}
\label{sec:local}
An alternative approach to the problem is given by traditional
(nonconvex) numerical optimization. There exist many variations, but
arguably the most successful methods for relatively small problems
such as the present ones are based on local gradient and Hessian
information. Typical algorithms in this class employ an iterative
scheme, combining the Newton search direction in combination with a
line search \cite{Nocedal}. These methods are reasonably fast in
converging to a \emph{local} minimum. For the larger problems in our
family, they usually converge to a stationary point within 10
seconds. However, they often end up in the wrong solution, unless a
very accurate starting point is given.

The drawbacks of local optimization methods are well-known: lacking
convexity, there are no guarantees of global (or even local)
optimality. Worse, even if in the course of the optimization we
actually reach the global minimum, there is usually no computationally
feasible way of verifying optimality.

Nevertheless, local optimization is an important tool for polynomial
problems, as is the use of homotopy methods to trace the optimal value
under small changes in the input data.  It would interesting to
investigate how these local numerical techniques can be best combined
with the computations to be described next.

\subsection{Experimental results using computational algebra}

In Table~\ref{tab:phctimes} we present typical running times for the
homotopy based approach, described in Section~\ref{sec:hom}. These
were obtained running PHCpack in ``black-box'' mode (\texttt{phc -b}),
that requires no user-specified parameters. The software traces all
solutions (not necessarily real), its number being equal to the
B\'ezout number.
Comparing with Table~\ref{tab:matsiz}, we can notice the adverse
effect of large B\'ezout numbers in the practical performance of the
algorithm, in spite of Verschelde's impressive implementation.

\begin{table}
\begin{center}
\begin{tabular}{c|rrrrrrrr}
$2d \,\,\, \backslash \,\,\, n$ & 3 & 5 & 7 & 9 & 11 & 13  & 15 \\ \hline
4 &  0.67 & 28.9  &  526 & - &  - &  - &  - \\
6 &  12.3 & 2643 &  - &  - &  - &  - &  - \\
8 &  70.6 & - &  - &  - &  - &  - &  - \\
10 &  508 & - &  - &  - &  - &  - &  - \\
\end{tabular}
\end{center}
\caption{Running time (in seconds) for the homotopy method.}
\label{tab:phctimes}
\end{table}

For the eigenvalue approach outlined in Section~\ref{sec:eigs}, we
compute the matrix ${\rm T}_f$ using a straightforward 
implementation in  Macaulay~2: the endomorphism ${\rm Times}_f$ is
constructed, and applied to the elements of the monomial
basis $\mathcal{B}$. The resulting matrix, in a sparse floating point
representation, is sent to a file for further processing.
We found that the construction of the matrix  ${\rm T}_f$ 
takes a surprisingly  long time.
for instance, it took Macaulay~2 over $10$ minutes to produce
the $125 \times 125$-matrix for $2d = 6$, $n = 3$.
The eigenvalue problem itself is solved using MATLAB;
it exploits the sparsity of the matrix, and runs reasonably fast.
However, it appears that even a more efficient implementation
of this method will not be able to compete with the timings
in Table \ref{tab:phctimes}, let alone the timings in
Table \ref{tab:sdptimes}.

After several discouraging attempts for small examples, we did
not pursue a full implementation for the
resultant-based methods sketched in Section 2.2.

\subsection{Experimental results using semidefinite programming}

We ran several instances of polynomials in the family described above,
for values of $K$ equal to 100, 1000, and 10000. In
Table~\ref{tab:sdptimes} the typical running times for the
semidefinite programming based approach on a single instance are
presented. These are fairly constant across instances, and no special
structure is exploited (besides what SeDuMi does internally).

The number of random instances for each combination of the parameters
is shown in Table~\ref{tab:numberofinstances}.  These values were
chosen in order to keep the total computation time for a given
category in the order of a few hours.

\begin{table}
\begin{center}
\begin{tabular}{c|ccccccc}
$2d \,\,\, \backslash \,\,\, n$ & 3 & 5 & 7 & 9 & 11 & 13  & 15 \\ \hline
4& 2000 &  2000 & 2000  &  200 & 20 &  20 &  2 \\
6& 2000 &  200 & 20 &  - &  - &  - &  - \\
8& 2000 &  20 & - &  - &  - &  - &  - \\
10& 2000 &  - & - &  - &  - &  - &  - 
\end{tabular}
\end{center}
\caption{Number of random instances in each category ($K=100, 1000, 10000$).}
\label{tab:numberofinstances}
\end{table}

Regarding the accuracy of the relaxation, in \emph{all cases tested}
the condition $f^{sos}=f^*$ was satisfied. As explained in the
previous section, this can be numerically verified by checking if the
solution of the corresponding SDP has rank one, from which a candidate
global minimizer is obtained. Evaluating the polynomial at this point
provides an upper bound on the optimal value, that can be compared
with the lower bound $f^{sos}$. In all our instances, the difference
between these two quantities was extremely small, and within the range
of numerical error.

As an additional check, when we used different methods for solving the
same instance, we have verified the solutions against each other. As
expected, the solutions were numerically close, in many cases up to
machine precision.

In particular, it is noted that the approach can handle in a
reasonable time (less than 35 min.) the case of a quartic polynomial
in thirteen variables. Our largest examples have the same degree
($2d=4$) and fifteen variables, correspond to an SDP with a matrix of
dimensions $136 \times 136$ with $3876$ auxiliary variables, and can
be solved in a few hours. A quick glance at the corresponding B\'ezout
number in Table~\ref{tab:matsiz} makes clear the advantages of the
presented approach.

\begin{table}
\begin{center}
\begin{tabular}{c|rrrrrrrr}
$2d \,\,\, \backslash \,\,\, n$ & 3 & 5 & 7 & 9 & 11 & 13  & 15 \\ \hline
4 &  0.2 & 0.5  &  4.4 & 52 &  361 &  1994 &  $27400^*$ \\
6 &  0.3 & 21.2 &  1046 &  - &  - &  - &  - \\
8 &  1.2 & 669 &  - &  - &  - &  - &  - \\
10 & 6.6 & - &  - &  - &  - &  - &  -
\end{tabular}
\end{center}
\caption{Running time (in seconds) for the semidefinite programs. The marked instance was solved on a different machine.}
\label{tab:sdptimes}
\end{table}

\section{What Next ?}

We have demonstrated that the sums of squares relaxation is
a powerful and practical technique in polynomial
optimization.  There are many open questions, both
algorithmic and mathematical, which are raised by 
our experimental results. One obvious question is
how often does it occur that $f^{sos} = f^*$~? This
can be studied for our simple model (\ref{eq:polyfam}),
or, perhaps better, for various natural probability 
measures on the space of polynomials of bounded degree.
This question is closely related
to  understanding the inclusion of the convex cone
of forms that are sums of squares inside the cone of
positive semidefinite forms. For the three-dimensional
family of symmetric sextics, this problem was studied 
in detail by  Choi, Lam and Reznick \cite{lam}.
Their work is an inspiration, but it also
provides a warning as to how difficult the general case will be,
even for  ternary sextics without symmetry.

We hope to pursue some of the following directions of inquiry 
in the near future.

\subsection{Sparseness and Symmetry}

Most polynomial systems one encounters
are sparse in the sense that there are only few monomials
with nonzero coefficients. Methods involving Newton polytopes,
such as sparse resultants \cite{gkzbook}  and polyhedral homotopies 
 \cite{ver}, are designed to deal with such problems.
Symmetry with respect to finite matrix groups is another
feature of many polynomial problems arising in practise.
The book of Gatermann \cite{gat} is an excellent first reference.

We wish to adapt our semidefinite programming approach to
input polynomials $f$ which are sparse or symmetric or both. 
For instance,  our polynomial example (\ref{RunningEx}) is 
both sparse and invariant under permutation of the variables $x,y,z$.
Both Newton polytope techniques and representation theory
can be used to reduce the size of the matrices 
and the number of free parameters in the semi-definite programs.

\subsection{Higher degree relaxations}

If we are unlucky, then the output produced by SeDuMi
will not satisfy the hypothesis of Proposition
\ref{completesquare}, and we conclude that the 
bound $f^{sos}$ is probably
strictly smaller than the optimal solution $f^*$. 
In that event we redo our computation in higher degree,
now with a  larger SDP.
The key idea is that even though $f(x)-\lambda$ may not be a
sum of squares, if there exists a positive  polynomial $g(x)$
such that $g (x) \cdot (f(x)-\lambda)$ is a sum of squares, then $\lambda \leq
f^*$. The choice of $g$ can be either made \emph{a priori} (for
instance, $g = \sum_{i=1}^n x_i^{2k}$), or as a result of an optimization
step (see \cite{pablo} for details). 
The Positivstellensatz (see Section 7) ensures that
$f^*$ will be found if the degree of $g$ is large enough.

\subsection{Solving polynomial equations} A natural application of 
Shor's relaxation, hinted at in \cite{shor},  is solving
polynomial systems $\, g_1(x) = \cdots = g_r(x) = 0$.
The polynomial $\,f(x) := \sum_{i=1}^r g^2_i(x)\,$ 
satisfies $\,f^* \geq f^{sos} \geq 0 $, and $f^* = 0$ holds if 
and only if the system has a real root.
Clearly,  $f^{sos} > 0$ is a sufficient condition for
the nonexistence of real roots.
An important open problem, essentially raised 
in \cite{shor}, is to characterize inconsistent systems  $\{g_1,\ldots,g_r\}$
with $f^{sos} = 0$. On the other hand, if $\,f^* = f^{sos} = 0 \,$ holds
then it is possible, at least in principle,
to obtain a numerical approximation of real roots using SDP.
However, for a robust implementation, perturbation arguments are
required and some important numerical issues arise, so the
perspectives for practical applications are still unclear.

\subsection{Minimizing polynomials over polytopes}
Consider a compact set
$$ P \quad = \quad \{ \, x \in \R^n \,\, : \,\,
\ell_1(x) \geq 0 , \ldots, 
\ell_s(x) \geq 0 \, \} ,$$
where $\ell_i$ is a linear form plus a constant,
say, $P$ is a polytope with $ s$ facets.
Handelman's Theorem \cite{han} states that
every polynomial  which is strictly positive on $P$
can be expressed as a positive linear combination 
of products $\, \ell_1(x)^{i_1} \cdots  \ell_s(x)^{i_s}$.
Suppose we wish to minimize a given polynomial function $f(x)$
over $P$. For $D \in \N$
we define the {\it $D$-th Handelman bound} $\,f^{(D)}\,$ to be
the largest $\, \lambda \in \R \,$ such that
$$ f(x) - \lambda \quad = 
\, \sum_{i_1 + \cdots + i_s \leq D} \!\!
 c_{i_1 \cdots i_s}  \cdot \ell_1(x)^{i_1} \cdots  \ell_s(x)^{i_s}
\qquad \hbox{for some} \,\,\,\,
 c_{i_1 \cdots i_s} \geq 0. $$
Handelman's Theorem states that the increasing sequence
$\,  f^{(D)},\, f^{(D+1)} ,\, f^{(D+2)}, \ldots \,$
converges to the  minimum of $f$ over $P$. Each bound
$f^{(D)}$ can be computed using {\bf linear programming only}.
It would be interesting to study the quality of these bounds,
and the running time of these linear programs,
and to see how things improve as we augment the approach
with semidefinite programming techniques.

\subsection{Which semialgebraic sets are semidefinite~?}
The feasible set of an SDP can be expressed by a linear 
matrix inequality, as in the dual formulation in Section 4.2.
It would be interesting to study these feasible sets 
using techniques from real algebraic geometry, and to
identify characteristic features of these sets.
Here is a very concrete problem whose solution, 
to the best of our knowledge, is not known. Fix three real symmetric 
matrices  $A,B$ and $C$  of size $N \times N$. Then
$$  \, S \quad = \quad
 \bigl\{ \, (x,y) \in \R^2 \,: \,
x A \, + \, y B + \, C \,\, \hbox{is positive
semidefinite} \,\bigr\} $$
is a closed, convex, semialgebraic subset of the plane $\R^2$.
The problem is to find a good characterization of those
subsets $S$. Given a semialgebraic subset $S \subset \R^2$
which is closed and convex, how to decide whether 
a ``semidefinite representation'' exists, and, in the affirmative case, 
how to find matrices $A,B,C$ of minimum size.

\section{Numerical Real Algebraic Geometry and The Positivstellensatz}

The first part of the above title refers
to a paper by Sommese and Wampler \cite{sommese}.
This paper and other more recent ones suggest that
numerical algorithms will play an increasingly
important role in computational (complex) algebraic geometry,
and that polynomial systems will become much more visible
in the context of {\sl Scientific Computation}. Along the
same lines, the
fastest software for computing Gr\"obner bases, due to 
Faug\'ere \cite{faugere}, no longer uses the Buchberger
algorithm but replaces it by sophisticated 
numerical linear algebra. Faug\'ere's scheme
$$
\hbox{\sl Gr\"obner Bases} \,\,\, + \,\,\, \hbox{\sl Numerical Linear Algebra}
\quad \rightarrow \quad \hbox{\sl Polynomial Problems
over $\C$}
$$
has the potential of entering the standard repertoire
of Scientific Computation.

Following \cite{pablo} we propose an analogous scheme
for the field of real numbers:
$$
\hbox{\sl Positivstellensatz} \,\,+ \,\,\hbox{\sl
Semidefinite Programming}
\,\,\, \rightarrow \,\,\, \hbox{\sl Polynomial Problems
over $\R$}. $$
In what follows, we shall explain this relationship
and why we see the Positivstellensatz as the main catalyst
for a future role of real algebra in scientific computation.

The Positivstellensatz \cite{BCR} is a common generalization
of Linear Programming Duality (for linear inequalities)
and  Hilbert's Nullstellensatz (for an algebraically closed field).
It states that, for a system of polynomial equations and inequalities,
either there exists  a solution in $\R^n$, or there exists a certain
polynomial identity which bears {\sl witness} to the fact that
no solution exists.  For instance, a single
polynomial inequality $f(x) < 0$ either has a solution $x \in \R^n$,
or there exists an identity $\,g(x) f(x) = h(x) \,$ where
$g$ and $h$ are sums of squares.  See \cite{boko} for an exposition of
the Positivstellensatz from the perspective of computational geometry.
Finding a witness by linear programming is proposed in \cite[\S 7.3]{boko}.

Here is our punchline, first stated in the dissertation
of the first author \cite{pablo}:
{\sl A Positivstellensatz witness of bounded degree can be
computed   by semidefinite programming.
Here we can also optimize linear parameters in the coefficients.}
This suggests the following algorithm for deciding
a system of polynomial equations and inequalities: decide whether 
there exists a witness for infeasibility of degree $\leq D$,
for some $D \gg 0$.
If our system is feasible, then we might like to minimize
a polynomial $f(x)$ over the solution set. The 
{\sl $D$-th SDP relaxation} would be to ask for 
the largest real number $\lambda$ such that
the given system together with the inequality
$\, f(x) - \lambda < 0 \,$ has an infeasibility
 witness of degree $D$. This generalizes what was proposed in
Sections 6.2, 6.3 and 6.4.

It is possible, at least in principle, to use 
an a priori bound for the degree  $D$ in the
Positivestellensatz, however, the currently known
bounds are still very large. Lombardi and Roy recently
announced a bound which is triply-exponential in
the number $n$ of variables. We hope that such bounds
can be further improved, at least for some natural
families of polynomial problems arising in optimization.

Here is a very simple example in the plane to illustrate 
our method:
\begin{equation}
f \,:=\, x-y^2+3 \,\, \geq \,\,0 \,, \qquad 
g \, := \, y+x^2+2 \, = \, 0.
\label{eq:expstz}
\end{equation}
By the Positivstellensatz, the system $\,\{f \geq 0, \, g = 0\}\,$
 has no solution  $\,( x,y) \in \R^2\,$ if and only if
there exist polynomials
 $s_1,s_2,s_3 \in \R[x,y]$ that satisfy the following:
\begin{equation}
s_1 + s_2 \cdot f + 1 + s_3 \cdot g \,\equiv \, 0 \, , \quad
\hbox{where $\,\, s_1$ and $s_2$ are sums of squares}.
\label{eq:pstz}
\end{equation}
The $D$-th SDP relaxation of the polynomial problem
$\,\{f \geq 0, \, g = 0\}\,$ asks whether there exists
a solution $(s_1,s_2,s_3)$ to
(\ref{eq:pstz}) where the polynomial $s_1$ has degree $\leq D$ and
the polynomials
$s_2,s_3$ have degree $\leq  D-2$. For each fixed integer $D > 0$
this can be tested by semidefinite programming.
For $D=2$ we find the  solution
\[
{\textstyle s_1 = \frac{1}{3} 
+ 2 \left( y+\frac{3}{2}\right)^2
+ 6 \left( x-\frac{1}{6}\right)^2}, \qquad
s_2 = 2, \qquad s_3 = -6.
\]
The resulting identity (\ref{eq:pstz})
proves the inconsistency of the system $\,\{f \geq 0, \, g = 0\}$.


\providecommand{\bysame}{\leavevmode\hbox to3em{\hrulefill}\thinspace}
\providecommand{\MR}{\relax\ifhmode\unskip\space\fi MR }
\providecommand{\MRhref}[2]{%
  \href{http://www.ams.org/mathscinet-getitem?mr=#1}{#2}
}
\providecommand{\href}[2]{#2}

\end{document}